\documentclass[amscd,amssymb,verbatim]{amsart}
\usepackage{amssymb,amsfonts,amsmath,amscd}

\theoremstyle{plain}
\newtheorem{Thm}{Theorem}[section]

\theoremstyle{definition}

\theoremstyle{remark}

\newcommand{\cC}{{\mathcal C}}

\newcommand{\N}{\mathbb{N}}

\newcommand{\R}{\mathbb{R}}

\newcommand{\bsa}{\boldsymbol{a}}
\newcommand{\bsb}{\boldsymbol{b}}
\newcommand{\bsc}{\boldsymbol{c}}
\newcommand{\bsd}{\boldsymbol{d}}

\newcommand{\bsf}{\boldsymbol{f}}
\newcommand{\bsF}{\boldsymbol{F}}
\newcommand{\bsG}{\boldsymbol{G}}

\newcommand{\bsH}{\boldsymbol{H}}

\newcommand{\bss}{\boldsymbol{s}}
\newcommand{\bst}{\boldsymbol{t}}

\newcommand{\bsx}{\boldsymbol{x}}

\newcommand{\bsy}{\boldsymbol{y}}

\newcommand{\bsz}{\boldsymbol{z}}
\newcommand{\bo}{{\mathbf 0}}

\begin{document}
\title[A General Implicit/Inverse Function Theorem]{A General Implicit/Inverse Function Theorem}
\author{Bruce Blackadar}
\address{Department of Mathematics/0084 \\ University of Nevada, Reno \\ Reno, NV 89557, USA}
\email{bruceb@unr.edu}

\date{\today}

\maketitle
\begin{abstract}
The Implicit and Inverse Function Theorems are special cases of a general Implicit/Inverse
Function Theorem which can be easily derived from either theorem.  The theorems can thus be
easily deduced from each other via the generalized version.
\end{abstract}

\section{Introduction}

The Implicit and Inverse Function Theorems are the deepest results from differential calculus, and among
the most important; they are the primary means of passing from the {\em infinitesimal} to the
{\em local} in analysis and geometry.  Although they are distinct results, each is a relatively easy corollary
of the other, so one only needs to carry out the (rather difficult) proof of one of them
to obtain both.  
They can both be regarded as aspects of the general problem of solving
systems of functional equations; existence and uniqueness theorems for solutions to differential equations
fall into the same class of results.

These theorems have a long and complicated history, going back at least to {\sc Descartes} (although for a long time
it was not always appreciated that they had to be {\em proved}), and there are many versions.  
The $\cC^1$ Implicit Function Theorem (\ref{SmoothImpFnThm}) was first proved by {\sc Dini} in the 1870's,
although {\sc Dini's} contributions were not properly recognized for a long time since he published them
only in his lecture notes, which had a limited circulation.  For continuous and higher-order smooth versions, see \cite{GoursatImplicites}
and \cite{YoungImplicit}.  The book \cite{KrantzPImplicit} contains the most thorough treatment available of
the various forms of the Implicit and Inverse Function Theorems (although rather incredibly this book
does not contain a complete correct proof of the basic Implicit Function Theorem as stated in 
\ref{SmoothImpFnThm} or in most Real Analysis texts).

It is well known that the two theorems can be easily deduced from each other, and the arguments
can be found in a number of Real Analysis texts.  But it seems not to be well known that there is
actually a general theorem of which both the Implicit and Inverse Function Theorems are 
special cases (a version of this general theorem appeared in \cite{NijenhuisStrong}), and that the
general theorem can be easily proved from either standard theorem essentially by the usual
arguments for deriving the other theorem.  Thus, once one of the theorems is proved, it is most
natural to derive the other by stating and proving the general version.  The aim of this article is
to state the general version (Theorem \ref{GenImpInvFnThm}) and show how to obtain it from the 
usual ones. 

\section{Statements of the Theorems}

The idea of the theorems is that if one has $n$ equations in $m$ unknowns, then there is
``usually'' a unique solution if $n=m$, many solutions if $n<m$, and no solutions if $n>m$;
if $n<m$, it should be possible to solve for $n$ of the variables as functions of the others
(specifying values for $m-n$ of the variables yields a system of $n$ equations in the remaining $n$ variables).
By elementary Linear Algebra, these statements are literally true if the equations are linear unless there is
degeneracy in the system.  Since smooth functions are ``approximately linear,'' it might
also be possible (and in fact {\em is}) to do the same locally in the smooth case.

Vector notation is a convenient way to rephrase systems of equations.  An equation in
unknowns $x_1,\dots,x_m$ can, by moving everything to one side of the equation, be written
as $F(x_1,\dots,x_m)=0$, where $F:\R^m\to\R$ is a function.  A set of $n$ equations
$$F_1(x_1,\dots,x_m)=0$$
$$F_2(x_1,\dots,x_m)=0$$
$$\cdots$$
$$F_n(x_1,\dots,x_m)=0$$
can be written as $\bsF(x_1,\dots,x_m)=\bo$, or $\bsF(\bsx)=\bo$, where $\bsF$ is the function from $\R^m$ (or a subset) to $\R^n$
with coordinate functions $(F_1,\dots,F_n)$, and $\bsx=(x_1,\dots,x_m)$.

In the case $m>n$, it is convenient to change notation slightly by setting $p=m-n$.  We would
like to be able to solve the system $\bsF(\bsx)=\bo$ for $n$ of the variables in terms of the other $p$ variables.
(Here the relative sizes of $n$ and $p$ do not matter; we could have $n<p$, $n=p$, or $n>p$.)
We can in general choose any $n$ of the variables to solve for in terms of the other ones;
we will notationally choose the last $n$ and rename them $y_1,\dots,y_n$, i.e.\ our equation to solve
is of the form
$$\bsF(x_1,\dots,x_p,y_1,\dots,y_n)=\bo$$
or $\bsF(\bsx,\bsy)=\bo$, where $\bsx=(x_1,\dots,x_p)\in\R^p$ and $\bsy=(y_1,\dots,y_n)\in\R^n$.
A {\em solution} on a subset $W$ of $\R^p$ is a function $\bsf:W\to\R^n$ such that the pair $(\bsx,\bsy)$ for $\bsy=\bsf(\bsx)$
satisfies the equation for all $\bsx\in W$, i.e.\ $\bsF(\bsx,\bsf(\bsx))=\bo$ for all $\bsx\in W$.  We say the function $\bsf$
is {\em implicitly defined} by the equation $\bsF(\bsx,\bsy)=\bo$.  It is rarely possible to find an explicit formula
for an implicitly defined function in the nonlinear case.

Examples with $n=p=1$ from elementary calculus show that one can only expect
the equation $\bsF(\bsx,\bsy)=\bo$ to define $\bsy$ as a continuous function of $\bsx$ locally, i.e.\ in some neighborhood
of a point $\bsa$ where $(\bsa,\bsb)$ satisfies the equation, and then only around points where there is some nondegeneracy
condition on $\bsF$.  If $n=p=1$ and $F$ is smooth, the graph of the equation $F(x,y)=0$ is a curve in $\R^2$, and this curve is in general only
locally the graph of a smooth function, and not even necessarily around every point of the curve; points where the curve crosses
itself or has a vertical tangent (or has even more complicated local behavior) must be excluded.


In the statements of the theorems, we will use the following notation.  Let $\bsF$ be a function from an open set $U$ in $\R^{p+n}$ to $\R^n$.
Write the coordinates of $\R^{p+n}$ as $(x_1,\dots,x_p,y_1,\dots,y_n)$, or $(\bsx,\bsy)$ with $\bsx=(x_1,\dots,x_p)\in\R^p$ and
$\bsy=(y_1,\dots,y_n)\in\R^n$ as above, and let $F_1,\dots,F_n$ be the coordinate functions of $\bsF$.  So, in coordinates,
$$\bsF(\bsx,\bsy)=(F_1(\bsx,\bsy),\dots,F_n(\bsx,\bsy))$$
$$=(F_1(x_1,\dots,x_p,y_1,\dots,y_n),\dots,F_n(x_1,\dots,x_p,y_1,\dots,y_n)))\ .$$
Let $(\bsx,\bsy)\in U$, and suppose all partial derivatives of the form $\frac{\partial F_i}{\partial y_j}(\bsx,\bsy)$ exist.  Write
$$D_{\bsy}\bsF(\bsx,\bsy)=\left [ \begin{array}{ccc} {\frac{\partial F_1}{\partial y_1}(\bsx,\bsy)} & \cdots & {\frac{\partial F_1}{\partial y_n}(\bsx,\bsy)} \\
\vdots & \ddots & \vdots \\ {\frac{\partial F_n}{\partial y_1}(\bsx,\bsy)} & \cdots & {\frac{\partial F_n}{\partial y_n}(\bsx,\bsy)} \end{array} \right ]$$
and, if all partial derivatives of the form $\frac{\partial F_i}{\partial x_j}(\bsx,\bsy)$ exist, write
$$D_{\bsx}\bsF(\bsx,\bsy)=\left [ \begin{array}{ccc} {\frac{\partial F_1}{\partial x_1}(\bsx,\bsy)} & \cdots & {\frac{\partial F_1}{\partial x_p}(\bsx,\bsy)} \\
\vdots & \ddots & \vdots \\ {\frac{\partial F_n}{\partial x_1}(\bsx,\bsy)} & \cdots & {\frac{\partial F_n}{\partial x_p}(\bsx,\bsy)} \end{array} \right ]$$
(note that $D_{\bsy}\bsF(\bsx,\bsy)$ is $n\times n$ and $D_{\bsx}\bsF(\bsx,\bsy)$ is $n\times p$ if they are defined).  Then if $\bsF$ is differentiable
at $(\bsx,\bsy)$, we have that the $n\times(p+n)$ Jacobian matrix $D\bsF(\bsx,\bsy)$ partitions as
$$D\bsF(\bsx,\bsy)=\left [ \begin{array}{ccc} \  & \vdots & \  \\ D_{\bsx}\bsF(\bsx,\bsy) & \vdots & D_{\bsy}\bsF(\bsx,\bsy) \\ \  & \vdots & \  \end{array} \right ] \ .$$

\paragraph{}
\begin{Thm} {\sc [Implicit Function Theorem]}\label{SmoothImpFnThm}
Let $U$ be an open set in $\R^{p+n}$ and $(\bsa,\bsb)\in U$.  Let $\bsF:U\to\R^n$
with $\bsF(\bsa,\bsb)=\bo$.  Suppose $\bsF$ is $\cC^r$ on $U$ for some $r$, $1\leq r\leq\infty$.
If $D_{\bsy}\bsF(\bsa,\bsb)$ is invertible, i.e.\ $det(D_{\bsy}\bsF(\bsa,\bsb))\neq0$, then there are:
\begin{enumerate}
\item[(i)] An open set $V$ in $\R^{p+n}$ containing $(\bsa,\bsb)$, with $V\subseteq U$,
such that $D_{\bsy}\bsF(\bsx,\bsy)$ is invertible for all $(\bsx,\bsy)\in V$;
\item[(ii)] An open set $W$ in $\R^p$ containing $\bsa$;
\item[(iii)] A unique function $\bsf:W\to\R^n$ such that
$\bsf(\bsa)=\bsb$ and, for all $\bsx\in W$, $(\bsx,\bsf(\bsx))\in V$ and
$$\bsF(\bsx,\bsf(\bsx))=\bo\ .$$
\end{enumerate}
Additionally, $\bsf$ is $\cC^r$ on $W$,
and for any $\bsx\in W$, we have
$$D\bsf(\bsx)=-D_{\bsy}\bsF(\bsx,\bsf(\bsx))^{-1}D_{\bsx}\bsF(\bsx,\bsf(\bsx))$$
(in particular, $D_{\bsy}\bsF(\bsx,\bsf(\bsx))$ is invertible for all $\bsx\in W$).
\end{Thm}

\paragraph{}
\begin{Thm} {\sc [Inverse Function Theorem]}\label{InverseFnThm}\index{Inverse Function Theorem!multivariable}
Let $U$ be an open set in $\R^n$, $\bsb\in U$, and $\bsf:U\to\R^n$ a $\cC^r$ function for some $r$, $1\leq r\leq\infty$.  If $D\bsf(\bsb)$ is invertible,
then there is a neighborhood $V$ of $\bsb$, $V\subseteq U$, such that
\begin{enumerate}
\item[(i)]  $\bsf$ is one-to-one on $V$.
\item[(ii)]  $\bsf(V)$ is open in $\R^n$.
\item[(iii)]  $D\bsf(\bsx)$ is invertible for all $\bsx\in V$.
\item[(iv)]  $\bsf^{-1}$ is $\cC^r$ on $\bsf(V)$ and, for all $\bsy\in \bsf(V)$,
$$D(\bsf^{-1})(\bsy)=[D\bsf(\bsf^{-1}(\bsy))]^{-1}\ .$$
\end{enumerate}
\end{Thm}

We will not discuss proofs of the individual theorems here (there are several standard ones).  See \cite{BlackadarReal} (or many other Real 
Analysis texts) for a full discussion.

\section{Twin Theorems and a Generalization}

The Inverse Function Theorem and Implicit Function Theorem 
are twin theorems in the sense that each is easily derivable from the other; the arguments
are much easier than proving either theorem from scratch.  In fact, either argument can be
used to prove a common generalization:

\paragraph{}
\begin{Thm} {\sc [General Implicit/Inverse Function Theorem]}\label{GenImpInvFnThm}
Let $n\in\N$, $p\in\N\cup\{0\}$, $\bsa\in\R^p$, $\bsb,\bsc\in\R^n$, $U$ an open neighborhood
of $(\bsa,\bsb)$ in $\R^{p+n}$,
$\bsF:U\to\R^n$ a $\cC^r$ function ($1\leq r\leq\infty$) with $\bsF(\bsa,\bsb)=\bsc$ and
$D_{\bsy}\bsF(\bsa,\bsb)$ invertible.  Then there are:
\begin{enumerate}
\item[(i)]  An open neighborhood $V$ of $(\bsa,\bsb)$ in $\R^{p+n}$, contained in $U$,
such that $D_{\bsy}\bsF(\bsx,\bsy)$ is invertible for all $(\bsx,\bsy)\in V$;
\item[(ii)]  An open neighborhood $W$ of $(\bsa,\bsc)$ in $\R^{p+n}$;
\item[(iii)] A unique function $\bsG:W\to\R^n$ satisfying $\bsG(\bsa,\bsc)=\bsb$ and,
for all $(\bsx,\bsz)\in W$, $(\bsx,\bsG(\bsx,\bsz))\in V$ and 
$$\bsF(\bsx,\bsG(\bsx,\bsz))=\bsz\ .$$
\end{enumerate}
Additionally, $\bsG$ is $\cC^r$ on $W$, and for all $(\bsx,\bsz)\in W$,
$$D\bsG(\bsx,\bsz)=\left [ \begin{array}{ccc} \  & \vdots & \  \\ -D_{\bsy}\bsF(\bsx,\bsG(\bsx,\bsz))^{-1}D_{\bsx}\bsF(\bsx,\bsG(\bsx,\bsz)) & \vdots & D_{\bsy}\bsF(\bsx,\bsG(\bsx,\bsz))^{-1} \\ \  & \vdots & \  \end{array} \right ] \ .$$
\end{Thm}

\smallskip

\noindent
Note that the Inverse Function Theorem is the case $p=0$  (so the $\bsa$ is absent).  The
Implicit Function Theorem is a corollary by taking $\bsc=\bo$ and setting
$\bsf(\bsx)=\bsG(\bsx,\bo)$.  In fact, the Implicit Function Theorem says that for certain fixed
$\bsz$ near $\bsc$, there is a function $\bsf_{\bsz}$ implicitly defined near $\bsa$ by $\bsF(\bsx,\bsy)=\bsz$,
which is $\cC^r$ in $\bsx$; we have $\bsG(\bsx,\bsz)=\bsf_{\bsz}(\bsx)$, and the additional content
of the above theorem is that the $\bsf_{\bsz}$ are defined for all $\bsz$ near $\bsc$ and ``vary smoothly'' in $\bsz$.  Thus Theorem \ref{GenImpInvFnThm} is a ``parametrized version'' of both
the Implicit and Inverse Function Theorems.

\smallskip

The situation is symmetric in $\bsF$ and $\bsG$.  $\bsG$ satisfies the hypotheses of the theorem,
so there is a $\cC^r$ function $\bsH$ defined near $(\bsa,\bsb)$ for which
$\bsH(\bsa,\bsb)=\bsc$ and $\bsG(\bsx,\bsH(\bsx,\bsy))=\bsy$ for $(\bsx,\bsy)$ near $(\bsa,\bsb)$.
By the formula, $D\bsH(\bsx,\bsy)=D\bsF(\bsx,\bsy)$ for all $(\bsx,\bsy)$ near $(\bsa,\bsb)$,
and since $\bsH(\bsa,\bsb)=\bsF(\bsa,\bsb)$, we have $\bsH=\bsF$ on a neighborhood
of $(\bsa,\bsb)$, i.e.\ $\bsH=\bsF$ (perhaps restricted to a smaller neighborhood of $(\bsa,\bsb)$).
In particular, we have $\bsG(\bsx,\bsF(\bsx,\bsy))=\bsy$ for $(\bsx,\bsy)$ near $(\bsa,\bsb)$
(this can be verified directly).

\bigskip

We will give simple proofs of \ref{GenImpInvFnThm} from both the Inverse Function Theorem and
the Implicit Function Theorem.

\subsection*{Proof From Inverse Function Theorem}

Assume the Inverse Function Theorem \ref{InverseFnThm}.  We will deduce \ref{GenImpInvFnThm}.

Suppose $U\subseteq\R^{p+n}$ is open, $(\bsa,\bsb)\in U$, and $\bsF:U\to\R^n$ is a $\cC^r$ function with $\bsF(\bsa,\bsb)=\bsc$
and $D_{\bsy}F(\bsa,\bsb)$ invertible.  Define
$$\tilde \bsF (\bsx,\bsy)=(\bsx,\bsF(\bsx,\bsy))$$
for $(\bsx,\bsy)\in U$.  Then $\tilde \bsF$ is a $\cC^r$ function from $U$ to $\R^{p+n}$.  We have
$$D\tilde\bsF(\bsx,\bsy)=\left [ \begin{array}{cc} I_p & 0 \\ D_{\bsx}\bsF(\bsx,\bsy) & D_{\bsy}\bsF(\bsx,\bsy) \end{array} \right ]$$
for $(\bsx,\bsy)\in U$; thus $D\tilde\bsF(\bsa,\bsb)$ is invertible.  So by the Inverse Function Theorem there is a
neighborhood $V$ of $(\bsa,\bsb)$ on which $\tilde\bsF$ is one-to-one, $W=\tilde\bsF(V)$ is open, and the inverse function
$\Phi:W\to V$ is $\cC^r$.  Since $\tilde\bsF(\bsa,\bsb)=(\bsa,\bsc)$, we have $(\bsa,\bsc)\in W$ and $\Phi(\bsa,\bsc)=(\bsa,\bsb)$.

Write $\Phi(\bsx,\bsz)=(\bsH(\bsx,\bsz),\bsG(\bsx,\bsz))$, where $\bsH:W\to\R^p$ and $\bsG:W\to\R^n$.  Then $\bsH$ and $\bsG$ are $\cC^r$, and $\bsG(\bsa,\bsc)=\bsb$.  Since for $(\bsx,\bsz)\in W$ we have
$$(\bsx,\bsz)=\tilde\bsF(\Phi(\bsx,\bsz))=\tilde\bsF(\bsH(\bsx,\bsz),\bsG(\bsx,\bsz))=(\bsH(\bsx,\bsz),\bsF(\bsH(\bsx,\bsz),\bsG(\bsx,\bsz)))$$
we obtain that $\bsH(\bsx,\bsz)=\bsx$ and
$$\bsF(\bsH(\bsx,\bsz),\bsG(\bsx,\bsz))=\bsF(\bsx,\bsG(\bsx,\bsz))=\bsz$$
for $(\bsx,\bsz)\in W$.

Uniqueness of $\bsG$ follows from the fact that $\tilde\bsF$ is one-to-one: if $\tilde\bsG$ is a function from $W$ to $\R^n$ with
$(\bsx,\tilde\bsG(\bsx,\bsz))\in V$ and $\bsF(\bsx,\tilde\bsG(\bsx,\bsz))=\bsz$ for all $(\bsx,\bsz)\in W$, then for $(\bsx,\bsz)\in W$ we have
$$\tilde\bsF(\bsx,\tilde\bsG(\bsx,\bsz))=(\bsx,\bsz)=\tilde\bsF(\bsx,\bsG(\bsx,\bsz))$$
so $\tilde\bsG(\bsx,\bsz)=\bsG(\bsx,\bsz)$.

Thus \ref{GenImpInvFnThm} is proved (the derivative formula is a straightforward exercise using
the Chain Rule).

\subsection*{Proof From Implicit Function Theorem}

Now assume the Implicit Function Theorem \ref{SmoothImpFnThm} holds.  We will deduce 
\ref{GenImpInvFnThm}.

Suppose $U\subseteq\R^{p+n}$ is open, $(\bsa,\bsb)\in U$, and $\bsF:U\to\R^n$ is as in the statement
of \ref{GenImpInvFnThm}.  Set $\bsc=\bsF(\bsa,\bsb)$.  Write 
$\R^{p+2n}\cong\R^p\times\R^n\times\R^n$, and set
$$U'=\{(\bsx,\bsz,\bsy)\in\R^{p+2n}:(\bsx,\bsy)\in U\}\ .$$
Then $U'$ is open in $\R^{p+2n}$, and $(\bsa,\bsc,\bsb)\in U'$. Define $\tilde\bsF:U'\to\R^n$ by
$$\tilde\bsF(\bsx,\bsz,\bsy)=\bsF(\bsx,\bsy)-\bsz\ .$$
Then $\tilde\bsF$ is $\cC^r$ on $U'$, $\tilde\bsF(\bsa,\bsc,\bsb)=\bo$, and
$$D\tilde\bsF(\bsx,\bsz,\bsy)=\left [ \begin{array}{ccc} D_{\bsx}\bsf(\bsx,\bsy) & -I_n & D_{\bsy}\bsF(\bsx,\bsy) \end{array} \right ]$$
for $(\bsx,\bsz,\bsy)\in U'$, and in particular $D_{\bsy}\tilde\bsF(\bsa,\bsc,\bsb)=D_{\bsy}\bsF(\bsa,\bsb)$ is invertible.   

We can thus solve the equation $\tilde\bsF(\bsx,\bsz,\bsy)=\bo$
implicitly for $\bsy$ as a function of $(\bsx,\bsz)$, i.e.\ there are open neighborhoods $W$ of
$(\bsa,\bsc)$ and $V$ of $(\bsa,\bsb)$ and a unique function $\bsG:W\to\R^n$ such that
$\bsG(\bsx,\bsz)\in V$ and $\tilde\bsF(\bsx,\bsz,\bsG(\bsx,\bsz))=\bo$ for all $(\bsx,\bsz)\in W$,
and $\bsG$ is $\cC^r$.
But  $\tilde\bsF(\bsx,\bsz,\bsG(\bsx,\bsz))=\bo$ is equivalent to $\bsF(\bsx,\bsG(\bsx,\bsz))=\bsz$.

\smallskip

Thus \ref{GenImpInvFnThm} is proved (the derivative formula is again a straightforward exercise).

\section{Weakening the Differentiability}

Various versions of the Implicit and Inverse Function Theorems have been given over the years.
The theorems hold verbatim with Euclidean spaces replaced by Banach spaces.  For the Inverse
Function Theorem, the continuous differentiability
hypothesis on $\bsf$ can be relaxed to
\begin{enumerate}
\item[(i)]  strong differentiability at $\bsa$ with invertible derivative there (\cite{LeachNote}, \cite{KnightStrong})

\begin{center}
or
\end{center}

\item[(ii)]  differentiability in a neighborhood of $\bsa$ with invertible derivative everywhere in the neighborhood (\cite{SaintRaymondLocal}, implicit in  \cite{ChernavskiiFinite}; cf.\ \cite{Tao}).
\end{enumerate}

In each of these cases we get a corresponding Generalized Implicit/Inverse Function Theorem
by the same argument as in the proof of \ref{GenImpInvFnThm}, and thus an Implicit Function
Theorem:

\paragraph{}
\begin{Thm} (cf.\ \cite{NijenhuisStrong})
Let $n\in\N$, $p\in\N\cup\{0\}$, $\bsa\in\R^p$, $\bsb,\bsc\in\R^n$, $U$ an open neighborhood
of $(\bsa,\bsb)$ in $\R^{p+n}$,
$\bsF:U\to\R^n$ a function which is strongly differentiable at $(\bsa,\bsb)$, with $\bsF(\bsa,\bsb)=\bsc$ and
$D_{\bsy}\bsF(\bsa,\bsb)$ invertible.  Then there are:
\begin{enumerate}
\item[(i)]  An open neighborhood $V$ of $(\bsa,\bsb)$ in $\R^{p+n}$, contained in $U$;
\item[(ii)]  An open neighborhood $W$ of $(\bsa,\bsc)$ in $\R^{p+n}$;
\item[(iii)] A unique function $\bsG:W\to\R^n$ satisfying $\bsG(\bsa,\bsc)=\bsb$ and,
for all $(\bsx,\bsz)\in W$, $(\bsx,\bsG(\bsx,\bsz))\in V$ and 
$$\bsF(\bsx,\bsG(\bsx,\bsz))=\bsz\ .$$
\end{enumerate}
Additionally, $\bsG$ is strongly differentiable at $(\bsa,\bsc)$, and 
$$D\bsG(\bsa,\bsc)=\left [ \begin{array}{ccc} \  & \vdots & \  \\ -D_{\bsy}\bsF(\bsa,\bsb)^{-1}D_{\bsx}\bsF(\bsa,\bsb) & \vdots & D_{\bsy}\bsF(\bsa,\bsb)^{-1} \\ \  & \vdots & \  \end{array} \right ] \ .$$
\end{Thm}

\paragraph{}
\begin{Thm} 
Let $n\in\N$, $p\in\N\cup\{0\}$, $\bsa\in\R^p$, $\bsb,\bsc\in\R^n$, $U$ an open neighborhood
of $(\bsa,\bsb)$ in $\R^{p+n}$,
$\bsF:U\to\R^n$ an everywhere differentiable function with $\bsF(\bsa,\bsb)=\bsc$ and
$D_{\bsy}\bsF(\bsx,\bsy)$ invertible for all $(\bsx,\bsy)\in U$.  Then there are:
\begin{enumerate}
\item[(i)]  An open neighborhood $V$ of $(\bsa,\bsb)$ in $\R^{p+n}$, contained in $U$;
\item[(ii)]  An open neighborhood $W$ of $(\bsa,\bsc)$ in $\R^{p+n}$;
\item[(iii)] A unique function $\bsG:W\to\R^n$ satisfying $\bsG(\bsa,\bsc)=\bsb$ and,
for all $(\bsx,\bsz)\in W$, $(\bsx,\bsG(\bsx,\bsz))\in V$ and 
$$\bsF(\bsx,\bsG(\bsx,\bsz))=\bsz\ .$$
\end{enumerate}
Additionally, $\bsG$ is everywhere differentiable on $W$, and for all $(\bsx,\bsz)\in W$,
$$D\bsG(\bsx,\bsz)=\left [ \begin{array}{ccc} \  & \vdots & \  \\ -D_{\bsy}\bsF(\bsx,\bsG(\bsx,\bsz))^{-1}D_{\bsx}\bsF(\bsx,\bsG(\bsx,\bsz)) & \vdots & D_{\bsy}\bsF(\bsx,\bsG(\bsx,\bsz))^{-1} \\ \  & \vdots & \  \end{array} \right ] \ .$$
\end{Thm}

Note that these theorems are technically not generalizations of \ref{GenImpInvFnThm} since both
the hypotheses and conclusion are weaker.

\section{Continuous and Mixed Versions}

There is a continuous version of the usual Implicit Function Theorem:

\paragraph{}
\begin{Thm} {\sc [Implicit Function Theorem, Continuous Version]}\label{ContImpFnThm}
Let $U$ be an open set in $\R^{p+n}$ and $(\bsa,\bsb)\in U$.  Let $\bsF:U\to\R^n$
with $\bsF(\bsa,\bsb)=\bo$.  Suppose $\bsF$ is continuous on $U$ and that $D_{\bsy}\bsF$ exists and is continuous on $U$,
i.e.\ that $\frac{\partial \bsF}{\partial y_j}$ exists and is continuous on $U$ for $1\leq j\leq n$.
If $D_{\bsy}\bsF(\bsa,\bsb)$ is invertible, i.e.\ $det(D_{\bsy}\bsF(\bsa,\bsb))\neq0$, then there are:
\begin{enumerate}
\item[(i)] An open set $V$ in $\R^{p+n}$ containing $(\bsa,\bsb)$, with $V\subseteq U$,
such that $D_{\bsy}\bsF(\bsx,\bsy)$ is invertible for all $(\bsx,\bsy)\in V$;
\item[(ii)] An open set $W$ in $\R^p$ containing $\bsa$;
\item[(iii)] A unique function $\bsf:W\to\R^n$ such that
$\bsf(\bsa)=\bsb$ and, for all $\bsx\in W$, $(\bsx,\bsf(\bsx))\in V$ and
$$\bsF(\bsx,\bsf(\bsx))=\bo\ .$$
\end{enumerate}
Additionally, $\bsf$ is continuous on $W$.
\end{Thm}

There is also a mixed version where $\bsF$ is differentiable with respect to some, but not necessarily
all, of the $x$-variables as well as all the $y$-variables.  We will notationally divide the $x$-variables into
$q$ $s$-variables and $m$ $t$-variables, where $q+m=p$, where the $t$-variables are the ones for
which $\bsF$ is differentiable.  We also say a function of $\bsy$
and possibly some additional variables is $\cC^r$ ($r\geq 1$) in $\bsy$ if all partial derivatives in
the $y$-variables of order $\leq r$ exist and are continuous.

\paragraph{}
\begin{Thm} {\sc [Implicit Function Theorem, Mixed Version]}\label{MixImpFnThm}
Let $U$ be an open set in $\R^{q+m+n}$ and $(\bsa,\bsd,\bsb)\in U$.  Let $\bsF:U\to\R^n$
with $\bsF(\bsa,\bsd,\bsb)=\bo$.  Suppose $\bsF$ is continuous on $U$ and $\cC^r$ on $U$ 
in $(\bst,\bsy)$ for some $r$, $1\leq r\leq\infty$.
If $D_{\bsy}\bsF(\bsa,\bsd,\bsb)$ is invertible, i.e.\ $det(D_{\bsy}\bsF(\bsa,\bsd,\bsb))\neq0$, then there are:
\begin{enumerate}
\item[(i)] An open set $V$ in $\R^{q+m+n}$ containing $(\bsa,\bsd,\bsb)$, with $V\subseteq U$,
such that $D_{\bsy}\bsF(\bss,\bst,\bsy)$ is invertible for all $(\bss,\bst,\bsy)\in V$;
\item[(ii)] An open set $W$ in $\R^{q+m}$ containing $(\bsa,\bsd)$;
\item[(iii)] A unique function $\bsf:W\to\R^n$ such that
$\bsf(\bsa,\bsd)=\bsb$ and, for all $(\bss,\bst)\in W$, $(\bss,\bst,\bsf(\bss,\bst))\in V$ and
$$\bsF(\bss,\bst,\bsf(\bss,\bst))=\bo\ .$$
\end{enumerate}
Additionally, $\bsf$ is $\cC^r$ in $\bst$ on $W$,
and for any $(\bss,\bst)\in W$, we have
$$D_{\bst}\bsf(\bss,\bst)=-D_{\bsy}\bsF(\bss,\bst,\bsf(\bss,\bst))^{-1}D_{\bst}\bsF(\bss,\bst,\bsf(\bss,\bst))$$
(in particular, $D_{\bsy}\bsF(\bss,\bst,\bsf(\bss,\bst))$ is invertible for all $(\bss,\bst)\in W$).
\end{Thm}

The usual Implicit Function Theorem (\ref{SmoothImpFnThm}) is the case $q=0$ and the continuous
version (\ref{ContImpFnThm}) is the case $m=0$.  The continuous and mixed versions of the
Implicit Function Theorem cannot be proved directly from the Inverse Function Theorem, and thus
the mixed version is more general.  Some of the usual proofs of the Implicit Function Theorem
such as the proof by induction on $n$ (which was historically the first proof) work essentially
verbatim to give the mixed version (cf.\ \cite{BlackadarReal}).  We then get a general Mixed Implicit/Inverse Function Theorem
by the same argument as before:

\paragraph{}
\begin{Thm} {\sc [General Implicit/Inverse Function Theorem, Mixed Version]}\label{GenMixImpInvFnThm}
Let $n\in\N$, $q,m\in\N\cup\{0\}$, $\bsa\in\R^q$, $\bsd\in\R^m$, $\bsb,\bsc\in\R^n$, $U$ an open neighborhood
of $(\bsa,\bsd,\bsb)$ in $\R^{q+m+n}$,
$\bsF:U\to\R^n$ a continuous function which is $\cC^r$ ($1\leq r\leq\infty$) in $(\bst,\bsy)$, with $\bsF(\bsa,\bsd,\bsb)=\bsc$ and
$D_{\bsy}\bsF(\bsa,\bsd,\bsb)$ invertible.  Then there are:
\begin{enumerate}
\item[(i)]  An open neighborhood $V$ of $(\bsa,\bsd,\bsb)$ in $\R^{q+m+n}$, contained in $U$,
such that $D_{\bsy}\bsF(\bss,\bst,\bsy)$ is invertible for all $(\bss,\bst,\bsy)\in V$;
\item[(ii)]  An open neighborhood $W$ of $(\bsa,\bsd,\bsc)$ in $\R^{q+m+n}$;
\item[(iii)] A unique function $\bsG:W\to\R^n$ satisfying $\bsG(\bsa,\bsd,\bsc)=\bsb$ and,
for all $(\bss,\bst,\bsz)\in W$, $(\bss,\bst,\bsG(\bss,\bst,\bsz))\in V$ and 
$$\bsF(\bss,\bst,\bsG(\bss,\bst,\bsz))=\bsz\ .$$
\end{enumerate}
Additionally, $\bsG$ is $\cC^r$ in $(\bst,\bsz)$ on $W$, and for all $(\bss,\bst,\bsz)\in W$,
$$D_{(\bst,\bsz)}\bsG(\bss,\bst,\bsz)=\left [ \begin{array}{ccc} \  & \vdots & \  \\ -D_{\bsy}\bsF(\bss,\bst,\bsG(\bss,\bst,\bsz))^{-1}D_{\bst}\bsF(\bss,\bst,\bsG(\bss,\bst,\bsz)) & \vdots & D_{\bsy}\bsF(\bss,\bst,\bsG(\bss,\bst,\bsz))^{-1} \\ \  & \vdots & \  \end{array} \right ] \ .$$
\end{Thm}

Theorem \ref{GenImpInvFnThm} is the case $q=0$, the Inverse Function Theorem is the case
$q=m=0$, and there is a continuous version with $m=0$.

It is unclear whether there are continuous or mixed versions of the theorems of Section 4.

\bibliography{impinvref}

\begin{thebibliography}{Gou03}

\bibitem[Bla]{BlackadarReal}
Bruce Blackadar.
\newblock {\em Real Analysis}.
\newblock http://wolfweb.unr.edu/homepage/bruceb/Meas.pdf.

\bibitem[{\v{C}}er64]{ChernavskiiFinite}
A.~V. {\v{C}}ernavski{\u\i}.
\newblock Finite-to-one open mappings of manifolds.
\newblock {\em Mat. Sb. (N.S.)}, 65 (107):357--369, 1964.

\bibitem[Gou03]{GoursatImplicites}
Ed. Goursat.
\newblock Sur la th\'eorie des fonctions implicites.
\newblock {\em Bull. Soc. Math. France}, 31:184--192, 1903.

\bibitem[Kni88]{KnightStrong}
William~J. Knight.
\newblock The {T}eaching of {M}athematics: {A} {S}trong {I}nverse {F}unction
  {T}heorem.
\newblock {\em Amer. Math. Monthly}, 95(7):648--651, 1988.

\bibitem[KP02]{KrantzPImplicit}
Steven~G. Krantz and Harold~R. Parks.
\newblock {\em The implicit function theorem}.
\newblock Birkh\"auser Boston Inc., Boston, MA, 2002.
\newblock History, theory, and applications.

\bibitem[Lea61]{LeachNote}
E.~B. Leach.
\newblock A note on inverse function theorems.
\newblock {\em Proc. Amer. Math. Soc.}, 12:694--697, 1961.

\bibitem[Nij74]{NijenhuisStrong}
Albert Nijenhuis.
\newblock Strong derivatives and inverse mappings.
\newblock {\em Amer. Math. Monthly}, 81:969--980, 1974.

\bibitem[Ray02]{SaintRaymondLocal}
Jean~Saint Raymond.
\newblock Local inversion for differentiable functions and the {D}arboux
  property.
\newblock {\em Mathematika}, 49(1-2):141--158 (2004), 2002.

\bibitem[Tao]{Tao}
Terence Tao.
\newblock
  http://terrytao.wordpress.com/2011/09/12/the-inverse-function-theorem-for-everywhere-differentiable-maps/.

\bibitem[You09]{YoungImplicit}
W.~H. Young.
\newblock On {I}mplicit {F}unctions and {T}heir {D}ifferentials.
\newblock {\em Proc. London Math. Soc.}, S2-7(1):397, 1909.

\end{thebibliography}
\bibliographystyle{alpha}

\end{document}